\documentclass[]{sn-jnl}

\usepackage{graphicx}%
\usepackage{multirow}%
\usepackage{amsmath,amssymb,amsfonts}%
\usepackage{amsthm}%
\usepackage{mathrsfs}%
\usepackage[title]{appendix}%
\usepackage{xcolor}%
\usepackage{textcomp}%
\usepackage{manyfoot}%
\usepackage{booktabs}%
\usepackage{algorithm}%
\usepackage{algorithmicx}%
\usepackage{algpseudocode}%
\usepackage{listings}%

\newlength\aptextwidth
\definecolor{BrickRed}{rgb}{0.588,0.098,0.055}

\def\noblue#1{\ifmmode \text{#1}\else #1\fi}
\def\noclh#1{\ifmmode \text{#1}\else #1\fi}
\def\rem#1{}

\def\ttable#1. #2{\begin{table}[t]\tablehat{#1}{#2}}
\def\mtable#1. #2{\begin{table}[hbtp]\tablehat{#1}{#2}}
\def\ptable#1. #2{\begin{table}[p]\tablehat{#1}{#2}}
\def\tablehat#1#2{\centering\small \vbox{\parindent=0pt
  \leftskip=0pt plus.5\hsize \rightskip=\leftskip \parfillskip=0pt
  ÒÀÁËÈÖÀ #1\\ #2}\nobreak\medskip\medskip }
\def\texendtable{\end{table}}
\def\notation{\par\ifnum\lastpenalty<25000 \bigbreak \fi
  \noindent\triangle\enspace\ignorespaces}

\def\hline{\noalign{\hrule}}

\let\bl=\bigl \let\br=\bigr

\def\eqnitemskip{\ifhmode \else \par
  \ifnum\lastpenalty>24999
    \ifnum\lastpenalty=25004 \fi
  \else \medbreak \fi \fi }
\def\eqnitem #1. {\eqnitemskip
  {\setbox0=\hbox{$#1^\circ$.\enspace}%
  \ifdim\wd0>\parindent \box0\ignorespaces \else
  \hbox to\parindent{\unhbox0\hss}\ignorespaces\fi}}
\def\eqnitemnobreak #1. {\noindent
  {\setbox0=\hbox{$#1^\circ$.\enspace}%
  \ifdim\wd0>\parindent \box0\ignorespaces \else
  \hbox to\parindent{\unhbox0\hss}\ignorespaces\fi}}
\newdimen\eqnparindent
\eqnparindent=\parindent
\def\eqnitem #1. {\eqnitemskip\noindent\hskip\eqnparindent $#1^\circ$.\enspace\ignorespaces }
\def\eqnitemnobreak #1. {\noindent\hskip\eqnparindent $#1^\circ$.\enspace\ignorespaces }
\def\simpleitem #1. {\eqnitemskip\noindent\hskip\eqnparindent #1.\enspace\ignorespaces }

\def\eqalignno#1{\displ@y \tabskip\centering
  \halign to\displaywidth{\hfil$\@lign\displaystyle{##}$\tabskip\z@skip
    &$\@lign\displaystyle{{}##}$\hfil\tabskip\centering
    &\llap{$\@lign\eqnofont##$}\tabskip\z@skip\crcr
    #1\crcr}}
\let\eqalignno=\eqalignm
\def\eqcenter#1{\displ@y \tabskip\centering
  \halign{\hfil$\displaystyle{##}$\hfil\crcr
    #1\crcr}}
\def\eqcenterno#1{\displ@y \tabskip\centering
  \halign to\displaywidth{\hfil$\@lign\displaystyle{##}$\hfil
    \tabskip\centering&\llap{$\@lign\eqnofont##$}\tabskip\z@skip\crcr
    #1\crcr}}
\def\texcases#1{\left\{\,\vcenter{\normalbaselines\m@th
    \ialign{$##\hfil$&\quad##\hfil\crcr#1\crcr}}\right.}
\def\Displaylines#1{\vcenter{\displ@y \tabskip\z@skip
  \halign{\hbox to\displaywidth{$\@lign\hfil\displaystyle##\hfil$}\crcr
    #1\crcr}}}

\def\tan{\mathop{\operator@font tg}\nolimits}

\begin{document}

\title[Notes on symmetries and reductions of algebraic equations]{Notes on symmetries and reductions of\\ algebraic equations}


\author*[1]{\fnm{Inna K.} \sur{Shingareva}}\email{inna.shingareva@unison.mx}

\author[2]{\fnm{Andrei D.} \sur{Polyanin}}\email{polyanin@ipmnet.ru}

\affil*[1]{\orgdiv{Department of Mathematics}, \orgname{University of Sonora},
\orgaddress{\street{Blvd. Luis Encinas y Rosales s/n}, \city{Hermosillo}, \postcode{83000}, \state{Sonora}, \country{Mexico}}}

\affil[2]{\orgdiv{Ishlinsky Institute for Problems in Mechanics}, \orgname{Russian Academy of Sciences},
\orgaddress{\street{101 Vernadsky Avenue, bldg 1}, \city{Moscow}, \postcode{119526}, \country{Russia}}}


\abstract{
Symmetries and reductions of some algebraic equations are considered.
Transformations that preserve the form of several algebraic equations,
as well as transformations that reduce the degree of these equations, are described.
Illustrative examples are provided.
The obtained results and solutions can be used as test problems 
for numerical methods of solving algebraic equations.
}

\keywords{algebraic equations, polynomial equations, symmetries, reductions, invariants, test problems for numerical methods}




\maketitle

\section{Introduction}\label{sec1}

For centuries, the main problem of algebra was considered to be
the development of methods for solving algebraic equations.
Formulas for solving algebraic equations of the first and second degrees
have been known since ancient times.

In the 16th century, Italian mathematicians Scipione del Ferro,
Niccol\'o Fontana (a.k.a. Tartaglia), Gerolamo Cardano, and Ludovico Ferrari obtained formulas
for solving algebraic equations of the third and fourth degree.
Historical information about this can be found in
\cite{Turnbull1947,Waerden1985}
(see also \cite{kor2000,polman2007,yac2012,the2016},
where various representations of solutions to these
equations are given).

Between the 16th and early 19th centuries,
for algebraic equations of degree higher than the fourth,
various attempts have been made to represent
{\it solutions in radicals}, i.e. in the form of an expression
containing only the coefficients of an algebraic equation and
the operations of addition, subtraction, multiplication, division and root extraction.
However, these attempts were unsuccessful in the general case.

In the 19th century, Paolo Ruffini, Niels Hendrik Abel, and \'Evariste Galois established
that in the general case the solution of algebraic equations of degree
higher than the fourth cannot be expressed in radicals
\cite{Waerden1985,Struik1986} (see also \cite{kin1996,pol2024}).

Note that the solution of individual algebraic equations
of the fifth degree and higher can be expressed by radicals.
Some such equations can be found in the handbook~\cite{pol2024},
where in addition to solutions by radicals, solutions are also given
that are expressed in terms of special functions
(for algebraic equations of the fifth degree, see also \cite{kin1996}).

In this note we will consider several algebraic equations that have certain symmetry properties,
which allows us to significantly simplify such equations.

\section{The simplest symmetries and transformations that reduce the order of algebraic equations}\label{sec2}

\eqnitem 1.
By symmetries of a mathematical equation we mean transformations
that preserve the form of the equation under consideration.
Reduction is a way of transforming a mathematical equation
into a simpler or more convenient (from some point of view)
form for analysis and solution.

Below we consider some algebraic equations that preserve their form
under certain linear transformations.
Let's start with a well-known simple example.
\medskip

\textit{Example 1.}
Consider the algebraic equation
\begin{equation}
a_{2n}x^{2n} + a_{2n-2}x^{2n-2} + a_{2n-4}x^{2n-4} + \cdots + a_{4}x^{4} + a_{2}x^{2} + a_{0} = 0,
\label{Eq2}
\end{equation}
containing only even powers.
The biquadratic equation is a special case of the equation \eqref{Eq2} for $n = 2$.
The substitution
\begin{equation}
x = -\bar x
\label{Eq3}
\end{equation}
leads to exactly the same equation for $\bar x$.
The equation \eqref{Eq2} is said to be invariant under the transformation \eqref{Eq3}.
It follows that if the equation \eqref{Eq2} has a solution $x = x_1$,
then it also has another solution $x = -x_1$.

By squaring both sides of \eqref{Eq3}, we obtain the simplest algebraic function,
which preserve its form under transformation \eqref{Eq3}:
\begin{equation}
x^2 = \bar x^2.
\label{Eq4}
\end{equation}
This function is the invariant of transformation \eqref{Eq3}.
If we choose invariant \eqref{Eq4} as a new variable, $z = x^2$,
then equation \eqref{Eq2} of order $2n$
is transformed to an equation of order $n$:
\begin{equation}
a_{2n}x^{n} + a_{2n-2}x^{n-1} + a_{2n-4}x^{n-2} + \cdots + a_{4}x^{2} + a_{2}x + a_{0} = 0.
\label{Eq5}
\end{equation}
Thus, in this case, the transition from the original variable $x$
to the transformation invariant \eqref{Eq4}, $z = x^2$, allows us
to simplify the equation under consideration (reduce its order by half).

\eqnitem 2.
Let us now consider a more general than \eqref{Eq5} algebraic equation of even degree~$2n$:
\begin{equation}
P_{2n}(x)=0,
\label{2.7a}
\end{equation}
which may contain lower terms with odd degrees.

{\sl {\bf Statement 1}.
Let the left-hand side of equation \eqref{2.7a} for some $\lambda$
and any $x$ identically satisfy the relation
\begin{equation}
P_{2n}(x)=P_{2n}(\lambda-x).
\label{2.7aa}
\end{equation}
Then the original equation, using the substitution $z=(x-\tfrac12\lambda)^2$,
is reduced to a simpler algebraic equation of degree $n$.}

Statement 1 is easy to prove using the relations
\begin{equation}
\begin{aligned}
P_{2n}(x)&=\tfrac12[P_{2n}(x)+P_{2n}(\lambda-x)]\\
&=\tfrac12\bl[P_{2n}\bl(\tfrac12\lambda+y\br)+P_{2n}\bl(\tfrac12\lambda-y\br)\br], \quad \ y=x -\tfrac12\lambda.
\end{aligned}
\label{2.7c}
\end{equation}
It follows from \eqref{2.7c} that equation \eqref{2.7a},
written in terms of the new variable $y$, does not change when $y$
is replaced by $-y$. Therefore, its left-hand side can be represented
as a linear combination of only even powers~$y$.
Next, the substitution $z=y^2$ ($z$ is the invariant of the transformation $\bar y=-y$)
is made, which halves the order of the equation.
\medskip

\textit{Example 2.}
Consider the equation
\begin{equation}
(x+a)^{2n}+(x+b)^{2n}-c=0,
\label{2.7b}
\end{equation}
which is a special case of equation \eqref{2.7a}.

It is easy to check that the left-hand side of \eqref{2.7b} does not change
if $x$ is replaced by $\lambda-x$, where $\lambda=-a-b$.
Therefore, the original equation, after raising the expressions on its
left-hand side to the $2n$ power and substituting $z=[x+\frac12(a+b)]^2$,
is reduced to an equation of degree $n$.
\medskip

\textit{Example 3.}
More complex than \eqref{2.7b} algebraic  equations
\begin{align*}
&(x+a)^{2n}+(x+b)^{2n}+c(x+a)(x+b)=d,\\
&(x+a)(x+b)[(x+a)^{2n}+(x+b)^{2n}]=c,
\end{align*}
which are special cases of equation \eqref{2.7a}
under condition \eqref{2.7aa}, are also simplified
by substituting $z=[x+\frac12(a+b)]^2$.

\eqnitem 3.
Let us now consider the algebraic equation of odd degree \text{$2n+1$}:
\begin{equation}
P_{2n+1}(x)=0.
\label{2.7a*}
\end{equation}

{\sl {\bf Statement 2}.
Let the left-hand side of equation \eqref{2.7a*} for some $\lambda$
identically satisfy the relation $P_{2n+1}(x)=-P_{2n+1}(\lambda-x)$.
Then the original equation has the root $x=\frac12\lambda$
and its left-hand side can be represented as
\begin{equation}
P_{2n+1}(x)=(x-\tfrac12\lambda)Q_{2n}(x),
\label{2.7a**}
\end{equation}
where a polynomial of even degree $Q_{2n}(x)$ identically satisfies
the relation $Q_{2n}(x)=Q_{2n}(\lambda-x)$.
According to Statement 1, the equation $Q_{2n}(x)=0$ is reduced to a simpler
algebraic equation of degree $n$ using
the substitution $z=(x-\tfrac12\lambda)^2$.}
\medskip

\textit{Example 4.}
Consider the equation
\begin{equation}
(x+a)^{2n+1}+(x+b)^{2n+1}+c(2x+a+b)=0,
\label{2.7b*}
\end{equation}
which is a special case of the equation \eqref{2.7a*}.

It is easy to check that the sign of the left-hand side of \eqref{2.7b*}
will change to the opposite if $x$ is replaced by $\lambda-x$, where $\lambda=-a-b$.
Therefore, the original equation has the root $x=\frac12\lambda=-\frac12(a+b)$
and its left-hand side can be represented as \eqref{2.7a**}.

\section{Reciprocal and related algebraic equations}\label{sec3}

\eqnitem 1.
Consider a \textit{reciprocal algebraic equation} of even degree of the form
\begin{equation}
a_0x^{2n} + a_1x^{2n-1} + a_2x^{2n-2} + \cdots + a_2x^2 + a_1x + a_0 = 0 \quad (a_0 \ne 0).
\label{Eq6}
\end{equation}
The left-hand side of this equation is called {\it a reciprocal polynomial}
(the coefficients of this polynomial, equally distant from its beginning and end, are equal).

In equation~\eqref{Eq6}, we make the substitution
\begin{equation}
x = \frac{1}{\bar x}.
\label{Eq7}
\end{equation}
After multiplying by $\bar x^{2n}$, we obtain exactly the same equation.
It follows that if equation \eqref{Eq6} has a solution $x = b$,
then it also has another solution $x = 1/b$.

It is easy to check that the simplest invariant of
transformation \eqref{Eq7} has the form
\begin{equation}
z = x + \frac{1}{x} = \bar x + \frac{1}{\bar x}.
\label{Eq11}
\end{equation}

{\bf Statement 3} (see, for example, \cite{Leung1992}).
{\sl A reciprocal equation of even degree \eqref{Eq6}
allows degree reduction by using substitution \eqref{Eq11}.
As a result, for the new variable $z$
we obtain an algebraic equation of degree $n$.}

\eqnitem 2.
Generalized reciprocal equations
\begin{align}
a_{0}\beta^nx^{2n}&+a_1\beta^{n-1}x^{2n-1}+\cdots+a_{n-1}\beta x^{n+1}+ a_nx^n\notag\\
&+\gamma a_{n-1}x^{n-1}+\gamma^2a_{n-2}x^{n-2}+\cdots+\gamma^{n-1}a_1x+\gamma^na_0 =0,\label{Eq10}
\end{align}
where $a_0,\beta,\gamma\not=0$, do not change during transformation
$$
x=\frac \gamma{\beta \bar x}.
$$
The simplest invariant of this transformation has the form
\begin{align}
z=\beta x+\frac \gamma{x}=\beta \bar x+\frac \gamma{\bar x}.
\label{Eq12}
\end{align}

{\bf Statement 4}.
{\sl Generalized reciprocal equation \eqref{Eq10}
allows degree reduction by using substitution~\eqref{Eq12}.
As a result, for the new variable~$z$ we obtain an algebraic equation of degree $n$.}
\medskip

\textit{Example 5.}
It is not difficult to show that the algebraic equation
\begin{equation*}
(x^2+a)^n+bx^{n-m}(x^2+a)^m+cx^n=0\qquad (n>m)
\end{equation*}
is a special case of equation \eqref{Eq10} and,
by using substitution \eqref{Eq12} with $\beta=1$ and $\gamma=a$,
this equation is reduced to the simpler equation
$\,z^n+bz^m+c=0$.

\section{Two algebraic equations that can be simplified}\label{sec4}

Below are two algebraic equations whose solution reduces to solving simpler equations.

\eqnitem 1.
Let us consider an algebraic equation of degree $n^2$,
\begin{equation}
a(ax^n+b)^n-x+b=0.
\label{ee02}
\end{equation}

This equation can be represented as
\begin{equation}
P(P(x))-x=0,
\label{ee03}
\end{equation}
where $P(x)=ax^n+b$ and
$P(P(x))$ is the second iteration of the binomial $P(x)$.
It is not difficult to verify that the roots of the simpler equation $P(x)-x=0$,
which is written in expanded form as $\,ax^n-x+b=0$, are also the roots of equation \eqref{ee02}.
Therefore, the {left-hand} side of \eqref{ee02} allows polynomial factorization and
can be represented as
\begin{equation}
P(P(x))-x=(P(x)-x)Q(x),
\label{ee04}
\end{equation}
where $Q(x)$ is a polynomial of degree $n^2-n$.
Thus, the solution of the original equation \eqref{ee02} is reduced to solving
two simpler equations, which are determined by equating the multipliers
to zero in the left part {of} \eqref{ee04}.
\medskip

\textit{Example 6.}
Let's consider an equation,
\begin{equation}
(x^2+b)^2-x+b=0,
\label{ee02*}
\end{equation}
which is a special case of equation \eqref{ee02} with $a=1$ and $n=2$ and
a special case of equation \eqref{ee03} with $P(x)=x^2+b$.
The {left-hand} side of \eqref{ee02*} allows polynomial
factorization and can be represented as the product of two quadratic polynomials
$$
(x^2+b)^2-x+b=(x^2-x+b)(x^2+x+b+1).
$$
The roots of these polynomials determine the roots {of} the original fourth-degree equation.

\eqnitem 2.
Consider another algebraic equation of degree $n^2$,
\begin{equation}
[a(ax+b)^n+b]^n-x=0.
\label{ee05}
\end{equation}

This equation can also be represented as \eqref{ee03}, where
$P(x)=(ax+b)^n$.
The roots of the simpler equation $P(x)-x=0$, which is written in expanded form as
\text{$\,(ax+b)^n-x=0$,} are also the roots of equation \eqref{ee05}.
Therefore, the {left-hand} side of \eqref{ee05} allows polynomial factorization
and can be represented as a product of polynomials \eqref{ee04}.
Thus, the solution of the original equation \eqref{ee05} is reduced
to solving two simpler algebraic equations.



\begin{thebibliography}{00}
\medskip

\setlength{\itemsep}{0.5em}


\bibitem{Turnbull1947}
\textit{Turnbull H.W.} Theory of Equations. Edinburgh: Oliver and Boyd, 1947.

\bibitem{Waerden1985}
\textit{Van der Waerden B.L.} A History of Algebra: From Al-Khwarizmi to Emmy Noether.
Berlin: Springer, 1985.


\bibitem{kor2000}
\textit{Korn G.A., Korn T.M.} Mathematical Handbook for Scientists and Engineers. New York: Dover Publ., 2000.


\bibitem{polman2007}
Polyanin A.D., Manzhirov A.V.  Handbook of Mathematics for
Engineers and Scientists. Boca Raton--London: Chapman \& Hall/CRC Press, 2007.

\bibitem{yac2012}
Yacoub M.D., Fraidenraich G.
A solution to the quartic equation. The Mathematical Gazette 2012, Vol. 96, No. 536, pp. 271--275.

\bibitem{the2016}
Tehrani F.T., Leversha G.
A simple approach to solving cubic equations. The Mathematical Gazette, 2016, Vol. 100 , No. 548, pp. 225--232.



\bibitem{Struik1986}
\textit{Struik D.J.} (ed.) A Source Book in Mathematics: 1200--1800.
Princeton: Princeton University Press, 1986.

\bibitem{pol2024}
\textit{Polyanin A.D.} Handbook of Exact Solutions to Mathematical Equations (Chapter~1). Boca Raton: CRC Press, 2024.

\bibitem{kin1996}
\textit{King R.B.} Beyond the Quartic Equation. Boston: Birkh\"auser, 1996.


\bibitem{Leung1992}
\textit{Leung K.T., Mok I.A.C., Suen S.N.} Polynomials and  Equations, Hong Kong: Hong Kong University Press, 1992.

\end{thebibliography}



\end{document}